\documentclass[12 pt]{article}
\usepackage[english]{babel}
\usepackage[active]{srcltx}
\usepackage{amsmath}
\usepackage{amsfonts,amssymb}
\usepackage[mathcal]{eucal}
\newtheorem{theorem}{Theorem}
\newtheorem{corollary}{Corollary}
\newtheorem{lemma}{Lemma}
\newtheorem{proposition}{Proposition}
\newtheorem{definition}{Definition}

\newcommand{\rav}{\stackrel{\triangle}{=}}
\newcommand{\ravref}[1]{\stackrel{(\ref{#1})}{=}}
\newcommand{\leqref}[1]{\stackrel{(\ref{#1})}{\leq}}
\newcommand{\geqref}[1]{\stackrel{(\ref{#1})}{\geq}}
\newcommand{\rref}[1]{$(\ref{#1})$}
\newcommand{\mm}[1]{{\mathbb{#1}}}
\newcommand{\av}{v}
\newcommand{\bw}{w}
\newcommand{\cnm}{{\mathbb{V}}^{\,i,\flat}}
\newcommand{\cnp}{{\mathbb{V}}^{\,i,\sharp}}
\newcommand{\cyg}{{\mathbb{V}}}
\newcommand{\avg}{\mathcal{V}}
\newcommand{\bwg}{\mathcal{W}}
\newcommand{\ssp}[2]{\left[#1\right]_{#2}}
\newcommand{\sspk}[2]{\Big[#1\Big]_{#2}}
\newcommand{\sspl}[1]{\Big[#1}
\newcommand{\sspll}[1]{\bigg[#1}
\newcommand{\ssplll}[1]{\Bigg[#1}
\newcommand{\sspr}[2]{#1\Big]_{#2}}
\newcommand{\ssprr}[2]{#1\bigg]_{#2}}
\newcommand{\ssprrr}[2]{#1\Bigg]_{#2}}
\newcommand{\ct}[1]{{\mathcal{#1}}}
\newcommand{\epsi}{\varepsilon}
\newcommand{\doc}{\em{Proof.}}
\newcommand{\bo}{\hfill {$\Box$}}

\begin{document}
\title{Tauberian theorem for value functions\thanks{Krasovskii Institute of Mathematics and Mechanics,     Russian
Academy of Sciences, 16, S.Kovalevskaja St., 620990, Yekaterinburg, Russia; \ \
Institute of Mathematics and Computer Science, Ural Federal University, 4, Turgeneva St., 620083, Yekaterinburg, Russia}}

\author{Dmitry Khlopin\\
{\it khlopin@imm.uran.ru} 
}

\maketitle





%

%
%
\maketitle

\begin{abstract}
      For two-person dynamic zero-sum games (both  discrete and continuous settings),  we investigate the limit of  value functions of finite horizon games with long run average cost as the time horizon tends to infinity and the limit of value functions of $\lambda$-discounted games as the discount tends to zero.
  We prove that the Dynamic Programming Principle for value functions directly leads to the  Tauberian Theorem---that the  existence  of a uniform limit  of  the value functions for one of the families implies  that the other one also uniformly converges to the same limit. No assumptions  on strategies  are necessary.
  To this end, we consider a mapping  that takes each payoff to the corresponding value function and preserves the sub- and super- optimality principles (the Dynamic Programming Principle). With their aid, we obtain certain inequalities on asymptotics of sub- and super- solutions, which lead to the  Tauberian Theorem. In particular, we consider the  case of differential games without relying on the existence of the saddle point;  a very simple stochastic game model is also considered.


	
\end{abstract}

\section{Introduction}\
 Hardy proved (see, for example,  \cite[Sect. 6.8]{Hardy1949})
  that,
for a bounded continuous function $g$, the limit of  long run averages  and the  limit of discounted averages (Cesaro mean  and
Abel mean, respectively)
  $$\frac{1}{T} \int_{0}^T g(t)\,dt,\qquad
    \lambda \int_{0}^\infty e^{-\lambda t}g(t)\,dt$$
    coincide if there exists at least one of these limits.
 This result and its generalizations have many applications (see, for example, \cite{Bingham,Feller,Korevaar}).

   Let us consider the analogs of this Tauberian theorem  for asymptotics of optimal values  in game-theoretic problem
   statements.
What if we optimize the Abel mean and/or Cesaro mean and then consider the limit of the optimal values corresponding to them?
Such a limit value (as the discount tends to zero) was first considered in
  \cite{Kh_Blackwell} for a stochastic formulation. As proved in
 \cite{BK}, for a stochastic two-person game with a finite number of states and actions, optimal long-time~averages and optimal
 discounted averages share the common limit. For more details on  the limit value for Abel mean and/or Cesaro mean in other
 stochastic formulations, see \cite{JN,SorinBig,LV,Renault2011,Sorin,Vigeral,Ziliotto2016}.

 In the deterministic case, the question of existence of limit values arose in the control theory, time and again; one may at the
 very least note
 \cite{chl1991,colklie,Gaitsgori1986}.  In the ergodic case (more generally, in the nonexpansive-like case) such limits exist
 and, moreover, they are usually independent of  the initial state, which was demonstrated in
\cite{ArisawaLions,Barles2000,wKAM1998}; although the results were released roughly at the same time, the methods of obtaining
them were thoroughly different.
 For the latest results on existence of limits of such values (for the nonexpansive-like case), see also
\cite[Sect. 3.4]{BCQ2013},\cite{CQ2015,QG2013,LQR,QuinRenault}.

A Tauberian theorem  (the equality of limit values)  was proved  for discrete time systems in  \cite{Lehrer}.
 The same result was obtained under an additional assumption that one of those limits is a constant function for control problems
 \cite{arisconst} and for differential games \cite{AlvBarditrue}.
    Note that, even in simple cases of control problems, the limits may not be constant functions
    \cite{grune,QuinRenault}.
In the general (non-ergodic) case, a
Tauberian theorem for very general  dynamic systems  was first proved in paper \cite{barton}.
 Then,   in \cite{Khlopin3},  a Tauberian theorem was proposed for differential games  (under Isaacs's condition).
 Later, in \cite{Ziliotto}, a very general approach  to proving Tauberian theorems was proposed for games  with two players with
 opposite goals
  in  discrete setting. In particular, it implies the same result for recursive games \cite{LV}.
Note that, in addition to uniform and exponential payoff families, the Tauberian theorems can be formulated for   arbitrary
probability densities. The corresponding results  are known for discrete time systems \cite{1993,Renault2013} and for optimal control and games \cite{CT15,LQR,Ziliotto2016}.

 The cornerstone of papers \cite{Khlopin3,barton} is the construction of near-optimal strategies for one of  the averages  by
 pieces of  near-optimal strategies for another average. It requires some assumptions on players' strategies, in particular, the
 Dynamic Programming Principle; in addition, for games,  there must be the existence of saddle point (see \cite{Khlopin3,CT15},
 and unpublished work \cite{KhlopinArXiv}).
  Paper \cite{Ziliotto} exhibits a more subtle approach. In stochastic games, the value is a fixed point of the Shapley operator
  for the corresponding game. Parameterized (by the discount or the finite horizon) families of the corresponding Shapley
  operators were embedded into certain Lipschitz continuous families of nonexpansive operators, and the corresponding Tauberian
  theorem was proved for the fixed points of the latter operators.

   The main aim of this paper is to obtain the Tauberian theorem as a direct consequence of the Dynamic Programming Principle without any technical assumptions on strategies, payoff functions, or anything else.
    To this end, we introduce a mapping (called a game value map) that assigns to every payoff the corresponding value function. In
  the general case,  this map can be constructed only  by payoff functions corresponding  to the Abel mean and Cesaro mean as
  payoffs, no strategies required. Considering a properly chosen chain of payoffs, using the monotonicity of the  game value map
  and sub- and super- optimality principles, we obtain  one-sided  inequalities on asymptotics  leading to all Tauberian
  theorems.
  Since no additional assumptions are imposed on the players' strategies (compare with \cite{Khlopin3,KhlopinArXiv}), neither topological nor measurable structures are used (compare with \cite{Ziliotto2016}) to prove the Tauberian theorem itself, not even the existence of saddle point (as in \cite{Khlopin3,CT15,Ziliotto}) is required, although a reduction to the typical formalization will apparently require some of these.

We also apply  this general theorem to zero-sum dynamic games  in continuous and discrete settings and to differential games
without a saddle point.

 The structure of the paper is as follows. We start by formulating the Tauberian theorem for dynamic games in continuous setting
 (Theorem~\ref{normal1}) in Sect.~\ref{abstractgame}. Then, we consider the general statement: we define the concept of a game
 value map and formulate the Tauberian theorem for this map (Theorem~\ref{normal2}) and the one-side inequalities on asymptotics
 (Propositions~\ref{av_bw}--\ref{bw_av_}).
    Sect.~\ref{proof} contains the proofs of   Theorems~\ref{normal1} and \ref{normal2}.
     Sections~\ref{apply} and \ref{diff} are devoted to the Tauberian theorem  for games in discrete setting
     (Theorem~\ref{normal3}) and  for  differential games
     (Theorem~\ref{maintheoremdiff}), respectively.  Also, the Tauberian theorem  for a very simple stochastic game model
     (Corollary~\ref{normal11}) is shown
     in  Sect.~\ref{apply}.
      The unwieldy and cumbersome proofs of propositions are confined to Appendix.

 \section{A dynamic zero-sum game}
\label{abstractgame}

 {\bf Dynamic system.}\
 Set $\mm{R}_+\rav\mm{R}_{\geq 0}.$   Assume the following items are given:
 \begin{itemize}
   \item   a nonempty set  $\Omega$, the state space;
   \item   a nonempty set $\mm{K}$ of maps from $\mm{R}_+$ to $\Omega$;
   \item   a running cost $g:\Omega\mapsto [0,1]$; for each process $z\in \mm{K}$, the  map $t\mapsto g(z(t))$ is  assumed to
       be Borel measurable.
 \end{itemize}

 {\bf On  payoffs.}\
 Let us now define a time average $\av_T(z)$ and a discount average $\bw_\lambda(z)$ for each process $z\in \mm{K}$
 by the following rules:
 \begin{eqnarray*}
 \av_T(z)\rav\frac{1}{T}\int_{0}^T g(z(t))\,dt,\quad
 \bw_\lambda(z)\rav\lambda\int_{0}^\infty e^{-\lambda t} g(z(t))\,dt\qquad \forall T,\lambda>0,z\in\mm{K}.
 \end{eqnarray*}
   Note that the definitions are valid, and  the means lie within  $[0,1].$

 {\bf On  lower games.}\
  For all $\omega\in\Omega$, let there be given non-empty sets $\ct{L}(\omega)$ and $\ct{M}(\omega).$
Let, for all   $\omega\in\Omega$, each pair $(l,m)\in\ct{L}(\omega)\times\ct{M}(\omega)$ of players' rules
generate a unique process $z[\omega,l,m]\in \mm{K}$  such that $z[\omega,l,m](0)=\omega$.

  The lower game is conducted in the following way:
  for a given $\omega\in\Omega,$
  the first player shows $l\in \ct{L}(\omega)$, and
  then, the second player chooses $m\in \ct{M}(\omega)$.
  The value function of this game is
  \begin{eqnarray}
 \cyg[c](\omega)\rav\sup_{l\in \ct{L}(\omega)}\inf_{m\in \ct{M}(\omega)}c(z[\omega,l,m])\qquad \forall\omega\in\Omega.
 \label{295}
\end{eqnarray}
  For instance, for every $T,\lambda>0$,
   the payoffs $\av_T,\bw_\lambda$ generate the following value functions:
  \begin{eqnarray*}
 \avg_T(\omega)&\rav&\cyg[\av_T](\omega)=\sup_{l\in \ct{L}(\omega)}\inf_{m\in \ct{M}(\omega)}\frac{1}{T}\int_{0}^T
 g(z[\omega,l,m](t))\,dt\qquad \forall\omega\in\Omega,\\
  \bwg_\lambda(\omega)&\rav&\cyg[\bw_\lambda](\omega)=\sup_{l\in \ct{L}(\omega)}\inf_{m\in \ct{M}(\omega)}\int_{0}^\infty \lambda
  e^{-\lambda t} g(z[\omega,l,m](t))\,dt\qquad
 \forall\omega\in\Omega.
\end{eqnarray*}

     Let us say that  the payoff family $\av_T (T>0)$
enjoys {\it the Dynamic Programming Principle} iff, for all $T>0$, the value function
$\avg_{T}$ coincides with the value functions for the following  payoffs:
 $$\mm{K}\ni z\mapsto \frac{1}{T}\int_{0}^h g(z(t))\,dt+\frac{T-h}{T}\avg_{T-h}(z(h))\quad \forall h\in(0,T).$$

     Let us say that the payoff family $\bw_\lambda (\lambda>0)$
enjoys {\it  the Dynamic Programming Principle} iff,  for all $\lambda>0$, the value function $\bwg_\lambda$
 coincides with the value functions for the following  payoffs:
$$\mm{K}\ni z\mapsto
 \lambda\int_{0}^h e^{-\lambda t}g(z(t))\,dt+e^{-\lambda h}\bwg_{\lambda}(z(h)) \quad\forall h>0.$$

\begin{theorem}\label{normal1}
Assume that the payoff families $\av_T (T>0)$ and  $\bw_\lambda (\lambda>0)$
enjoy the Dynamic Programming Principle.

 Then, the following two statements are equivalent:
\begin{description}
  \item[$(\imath)$]   The family of functions $\avg_T$ $(T>0)$  converges uniformly on $\Omega$ as $T\uparrow\infty.$
  \item[$(\imath\imath)$] The family of functions $\bwg_\lambda$ $(\lambda>0)$ converges uniformly on $\Omega$ as $\lambda\downarrow0$.
\end{description}
Moreover, when at least one of these statements holds, we have
    \begin{equation*}
    \lim_{T\uparrow\infty}\avg_T(\omega)=
       \ \lim_{\lambda\downarrow 0}\bwg_\lambda(\omega)\quad \forall\omega\in\Omega.
       \end{equation*}
\end{theorem}
For the proof of this theorem, refer to Sect.~\ref{proof}.

{\bf On  abstract control systems.}\
We can obtain
the Tauberian theorem for an abstract control system.
Following \cite{barton}, assume the sets $\Omega$, $\mm{K}$ to be given; for all $\omega\in\Omega$, let $\ct{L}(\omega)$ be the
set of all  feasible processes $z\in\mm{K}$
that begin at $\omega$. Let
$\ct{M}(\omega)$ be a singleton for all $\omega\in\Omega$.

Now, Theorem~\ref{normal1} implies
   \begin{corollary}
     \label{u1}
Assume that the payoff families $\av_T (T>0)$ and  $\bw_\lambda (\lambda>0)$
enjoy the Dynamic Programming Principle.


 Then, the following two statements are equivalent:
\begin{description}
  \item[$(\imath)$]   The maps $\Omega\ni\omega\mapsto
\sup_{z\in\ct{L}(\omega)} \av_T(z)$
      converge uniformly on $\Omega$ as $T\uparrow\infty.$
  \item[$(\imath\imath)$] The maps $\Omega\ni\omega\mapsto\sup_{z\in\ct{L}(\omega)} \bw_\lambda(z)$
       converge uniformly on $\Omega$ as $\lambda\downarrow0$.
\end{description}
Moreover, when at least one of these statements holds, we have
    \begin{equation*}
    \lim_{T\uparrow\infty}\sup_{z\in\ct{L}(\omega)} \av_T(z)=
       \lim_{\lambda\downarrow 0}\sup_{z\in\ct{L}(\omega)} \bw_\lambda(z)\quad \forall\omega\in\Omega.
       \end{equation*}
    \end{corollary}
  In \cite{barton}, it is stated that the Tauberian theorem holds for an abstract control system  if $\mm{K}$
is closed with respect to concatenation. This condition can be refined,  see \cite{khlopin2016}.

Theorem~\ref{normal1} and Corollary~\ref{u1} look similar to  Tauberian theorems for dynamic zero-sum game and control system with continuous setting, respectively, in   the most general statement.
Nevertheless, let us sketch the examples when this similarity is misleading.  In certain  game problems with information discrimination, the second player has to choose~$m$ from $\ct{M}(\omega,l)$ instead of from $\ct{M}(\omega)$ \cite{KhCh2000}. Value functions with $\inf_{l}\sup_{m}\inf_{\tau}$ instead of
$\sup_{l}\inf_{m}$ are applied in Hamilton-Jacobi-Isaacs variational inequalities  \cite[(17.7)]{subb} and
pursuit-evasion-defense problems \cite[(12)]{Fisac}.  Finally,   for instance,  for the control problem,  a  maximization
 of the expectation of the payoff with respect to some probability distribution
 (see, for instance, \cite[(11.3.2)]{oksendal}, \cite[(2.5)]{lacker})  is not covered by Corollary~\ref{u1}. For  this reason, in the next section, we introduce a mapping (called a game value map) that assigns to each payoff  the corresponding value function, and then we  formulate the Tauberian theorem for game value map.

 \section{General statement.}
\label{side}

{\bf On game value map.}\
Let the sets $\Omega$ and $\mm{K}$, running cost $g$, and payoffs $\av_T,\bw_\lambda$ be as before.
Denote by $\mathfrak{U}$ the set of all bounded maps from $\Omega$ to $\mm{R}$;  denote  by $\mathfrak{C}$ a non-empty set of maps from $\mm{K}$ to $\mm{R}$.
Thereinafter, the set  $\mathfrak{C}$ incorporates all conceivable payoffs,
and the set  $\mathfrak{U}$ contains all value functions for all games with  payoffs $c\in\mathfrak{C}.$

\begin{subequations}
Let $\mathfrak{C}$  satisfy
the following condition:
\begin{equation}
 Ac+B\in\mathfrak{C}\ \mathrm{  for\ all }\ A\geq 0,B\in\mm{R}\ \mathrm{ if } \ c\in\mathfrak{C}.
 \label{conditions}
\end{equation}
 In this section, we also assume that $\av_T,\bw_\lambda\in\mathfrak{C}$ for all positive $\lambda,T.$

A map $V$ from $\mathfrak{C}$ to $\mathfrak{U}$ is called {\it a game value map} if the following conditions hold:
\begin{eqnarray}
\label{conditions1}
  &\ &V[Ac+B]=A\,V[c]+B\  \textrm{ for all } c\in\mathfrak{C}, A\geq 0, B\in\mm{R},\\
  &\ &V[c_1](\omega)\leq V[c_2](\omega)\ \textrm{ for all }  \omega\in\Omega\ \textrm{ if } c_1(z)\leq c_2(z) \ \textrm{for all }
  z\in\mm{K}.
\label{conditions2}
\end{eqnarray}
\end{subequations}

{\bf On Dynamic Programming Principle.}\
 For all positive $\lambda,T,h>0$ and every function $U_*:\Omega\to\mm{R}$, define payoffs $\zeta^{U_*}_{h,T}:\mm{K}\to\mm{R}$,
 $\xi^{U_*}_{h,\lambda}:\mm{K}\to\mm{R}$ as follows:
 \begin{eqnarray*}
 \zeta^{U_*}_{h,T}(z)&\rav&
 \frac{1}{T+h}\int_{0}^h g(z(t))\,dt+\frac{T}{T+h}{U_*}(z(h))\qquad\forall  z\in\mm{K};\\
\xi^{U_*}_{h,\lambda}(z)&\rav&
\lambda\int_{0}^h e^{-\lambda t}g(z(t))\,dt+e^{-\lambda h}{U_*}(z(h))\qquad\forall z\in\mm{K}.
\end{eqnarray*}

\begin{definition}
  For a  game value map $V$,  let us say that a family of $U_T\in\mathfrak{U} (T>0)$
  is
  {\it a  subsolution}  \big({\it a  supersolution}\big) for
    the family of payoffs $\av_T (T>0)$
 if,
 for every $\varepsilon>0$,
  there exists natural  $\bar{T}$ such that,
     for all natural $h,T>\bar{T}$, the payoff
$\zeta^{U_T}_{h,T}$ lies in $\mathfrak{C}$ and enjoys
  $$U_{T+h}\leq V[\zeta^{U_T}_{h,T}]+\epsi\ \
  \Big(U_{T+h}\geq V[\zeta^{U_T}_{h,T}]-\epsi\Big).$$

  For a  game value map $V$,  let us say that a family of $U_\lambda\in\mathfrak{U} (\lambda>0)$
  is
  {\it a  subsolution}  \big({\it a  supersolution}\big)  for
  the family of payoffs $\bw_\lambda (\lambda>0)$
 if, for every $\varepsilon>0$,
  there exists natural $\bar{T}$ such that,
     for all natural $h>\bar{T}$ and positive $\lambda<1/\bar{T}$, the payoff
                        $\xi^{U_\lambda}_{h,\lambda}\in\mathfrak{C}$  lies in $\mathfrak{C}$ and enjoys
     $$U_{\lambda}\leq V[\xi^{U_\lambda}_{h,\lambda}]+\epsi\ \
  \Big(U_{\lambda}\geq V[\xi^{U_\lambda}_{h,\lambda}]-\epsi\Big).$$
    \end{definition}
For similar definitions, refer to
the  suboptimality principle \cite[Definition III.2.31]{Bardi} (also referred to as `stability with respect to second
player'~\cite{ks}) and \cite[Sect.~VI.4]{Bardi} for discrete setting.

\begin{definition}
 For a  game value map $V$,  let us say that the family of payoffs $\av_T(T>0)$ (resp., $\bw_\lambda(\lambda>0)$)
 enjoys {\it the weak Dynamic Programming Principle}   iff their
 value functions ($V[\av_T]$ and, resp., $V[\bw_\lambda]$) are, at the same time, a  subsolution and a supersolution for this
 payoff family.
\end{definition}

In particular, the family of payoffs $\av_T(T>0)$ (resp., $\bw_\lambda(\lambda>0)$)
enjoys {the weak Dynamic Programming Principle} if
$$V[\av_{T+h}]=V[\zeta^{V[\av_T]}_{h,T}],\quad\Big(V[\bw_{\lambda}]=V[\xi^{V[\bw_\lambda]}_{h,\lambda}]\Big)\qquad
\forall h,T\in\mm{N}, \lambda>0.$$

\begin{theorem}\label{normal2}
Let there be given a game value map $V:\mathfrak{C}\to\mathfrak{U}.$ Let  $\av_T,\bw_\lambda\in\mathfrak{C}$ for all
$\lambda,T>0.$
Assume that payoffs $\av_T (T>0)$ and  payoffs $\bw_\lambda (\lambda>0)$
enjoy the weak Dynamic Programming Principle.

Then, the following two statements are equivalent:
\begin{description}
  \item[$(\imath)$]   The family of functions $V[\av_T]$ $(T>0)$  converges uniformly on $\Omega$ as $T\uparrow\infty.$
  \item[$(\imath\imath)$] The family of functions $V[\bw_\lambda]$ $(\lambda>0)$ converges uniformly on $\Omega$ as $\lambda\downarrow0$.
\end{description}
Moreover, when at least one of these statements holds, we have
    \begin{equation*}
    \lim_{T\uparrow\infty}V[\av_T](\omega)=\lim_{\lambda\downarrow 0}V[\bw_\lambda](\omega)\quad \forall\omega\in\Omega.
       \end{equation*}
\end{theorem}
This theorem generalizes Theorem~\ref{normal1}. For its proof, refer to Sect.~\ref{proof}.

{\bf One-sided Tauberian theorems for bounds from above.}\
    In Appendix \ref{A}, we prove  the following proposition:
 \begin{subequations}
     \begin{proposition}
     \label{av_bw}
       For a game value map $V$,
       let a  family of functions
       $U_T\in\mathfrak{U} (T>0)$ and the family of functions $V[\bw_\lambda](\lambda> 0)$ be
       a  subsolution for payoffs $\av_T$ and
       a  supersolution for payoffs $\bw_\lambda $, respectively.

      Further, let $U_T (T>0)$
      satisfy
 \begin{equation}
 \label{slowlyT}
 \limsup_{T\uparrow\infty}\sup_{p\in[1,p_0]}\sup_{\omega\in\Omega}\,
\big(\,U_{T}(\omega)-U_{Tp}(\omega)\big)\leq 0.\quad\forall p_0>1.
 \end{equation}

      Then,
      for every $\varepsilon>0$, there exists a natural $N$ such that
       \begin{equation*}
       V[\bw_{\lambda}](\omega)\geq U_{1/\lambda}(\omega)-\varepsilon\qquad\forall \omega\in\Omega,\lambda\in(0,1/N).
       \end{equation*}
\end{proposition}

       In Appendix \ref{B}, we prove  a similar proposition, where the payoff families are swapped:
        \begin{proposition}
 \label{bw_av}
       For a game value map $V$,
       let a family of functions
       $U_\lambda\in\mathfrak{U} (\lambda>0)$ and the family of functions $V[\av_T](T> 0)$ be
       a  subsolution for payoffs $\bw_\lambda$ and
       a  supersolution for payoffs $\av_T$, respectively.

      Further, let
       $U_\lambda (\lambda>0)$
 satisfy
  \begin{equation}
 \label{slowlyl}
  \limsup_{\lambda\downarrow 0}\sup_{p\in[1,p_0]}
 \sup_{\omega\in\Omega}\,\big(\,U_{\lambda}(\omega)-U_{p\lambda}(\omega)\big)\leq 0\qquad \forall p_0>1.
 \end{equation}

       Then, for every
       $\varepsilon>0$,  there exists a natural $N$ such that
       $$V[\av_T](\omega)\geq U_{1/T}(\omega)-\varepsilon\qquad  \forall \omega\in\Omega,T>N.$$
     \end{proposition}

The inequalities similar to \rref{slowlyT} will be found in other Tauberian Theorems,
(see \cite[Definition 4.1.4]{Bingham},\cite[Sect.6.2]{Hardy1949},\cite[Theorem 4.1]{MN1981}).
 Also, we may take $p=p_0$ instead of $\sup_{p\in[1,p_0]}$
 in \rref{slowlyT} and \rref{slowlyl} for
    $U_T\rav V[\av_T]$ and  $U_\lambda\rav V[\bw_\lambda]$ using \rref{2012} and \rref{2022}  (see Appendix~\ref{AA}),
    respectively.
    On the other hand, conditions
  \rref{slowlyT} and \rref{slowlyl} are  tight enough for
    $U_T\rav V[\av_T]$ and  $U_\lambda\rav V[\bw_\lambda]$, respectively. In particular, if in Proposition~\ref{bw_av} one
    replaces \rref{slowlyl} with the following condition
  $$\lim_{p_0\downarrow 1}\limsup_{\lambda\downarrow 0}\sup_{p\in[1,p_0]}
 \sup_{\omega\in\Omega}\,\big|\,V[\bw_{\lambda}](\omega)-V[\bw_{p\lambda}](\omega)\big|=0,
$$
which holds for all game value maps from \rref{2022},
then Proposition~\ref{bw_av} would fail for
 a certain game value map (see, for example, \cite[Sect. 4]{barton} for control problems).

{\bf One-sided Tauberian theorems for bounds from below.}\
Applying these  propositions to $\mathfrak{C}^*\rav \{-c\,|\,c\in\mathfrak{C}\}$, $V^*[c]\equiv -V[-c]$, $g^*\equiv 1-g$ instead of $\mathfrak{C},V,g$, we obtain:
     \begin{proposition}
     \label{av_bw_}
      For a game value map $V$, 
          let a  family of functions
       $U_T\in\mathfrak{U} (T>0)$
        and the family of functions $V[\bw_\lambda](\lambda> 0)$ be
       a  supersolution for the payoffs $\av_T$ and
       a  subsolution for the payoffs $\bw_\lambda $, respectively.

      Further, let $U_T (T>0)$
      satisfy
 \begin{equation}
 \label{slowlyT_}
 \limsup_{T\uparrow\infty}\sup_{p\in[1,p_0]}\sup_{\omega\in\Omega}\,
\big(\,U_{Tp}(\omega)-U_{T}(\omega)\big)\leq 0.\quad\forall p_0>1.
 \end{equation}

      Then,
      for every $\varepsilon>0$, there exists a natural $N$ such that
       \begin{equation*}
       V[\bw_{\lambda}](\omega)\leq U_{1/\lambda}(\omega)+\varepsilon\qquad\forall \omega\in\Omega,\lambda\in(0,1/N).
       \end{equation*}
\end{proposition}
 \begin{proposition}
 \label{bw_av_}
      For a game value map $V$,
       let a family of functions
       $U_\lambda\in\mathfrak{U} (\lambda>0)$ and the family of functions $V[\av_T](T> 0)$ be
       a  supersolution for the payoffs $\bw_\lambda$ and
       a  subsolution for the payoffs $\av_T $, respectively.

      Further, let
       $U_\lambda (\lambda>0)$
  satisfy
  \begin{equation}
 \label{slowlyl_}
  \limsup_{\lambda\downarrow 0}\sup_{p\in[1,p_0]}
 \sup_{\omega\in\Omega}\,\big(\,U_{p\lambda}(\omega)-U_{\lambda}(\omega)\big)\leq 0\qquad \forall p_0>1.
 \end{equation}

       Then, for every
       $\varepsilon>0$,  there exists a natural $N$ such that
       $$V[\av_T](\omega)\leq U_{1/T}(\omega)+\varepsilon\qquad  \forall \omega\in\Omega,T>N.$$
     \end{proposition}
 \end{subequations}

For simplicity we could define $\mathfrak{C}$ as the set of all bounded maps from $\mm{K}$ to $\mm{R}$.
 It would be sufficient for proofs of all theorems of this article. We do not do it due to the following reasons.
 First, in stochastic  frameworks, each payoff $c$ has to be measurable with respect to some measurable space,
 see, for instance, Corollary~\ref{normal11} in Sect.~\ref{apply}.
 Secondly, all proofs use merely the boundedness of $\av_T,\bw_\lambda$, and  the additional requirement does not appear to help to obtain the bounds.

 Also note that in the definitions of  subsolution  and  supersolution, we consider only natural $h$ and $T$.  However, this strengthening is useless for  continuous setting,  whereas in discrete-time setting it allows a direct usage of the corresponding Dynamic Programming Principle  (see Section~\ref{apply}).

\section{Proofs of main results}
\label{proof}
{\bf Proof of Theorem~\ref{normal2}.}\
Since the payoffs  $\av_T (T>0)$ and the payoffs  $\bw_\lambda (\lambda>0)$
enjoy the weak Dynamic Programming Principle, we  see that $V[\av_T] (T>0)$ and $V[\bw_\lambda](\lambda>0)$ are
 simultaneously super- and subsolutions with respect to $V$
for the payoffs $\av_T$ and for the payoffs  $\bw_\lambda$ respectively.

If at least one of the considered limits  exists and is uniform   in $\Omega,$
then,
    either the limit of $V[\av_T]$ as $T\uparrow\infty,$
or the limit of $V[\bw_\lambda]$ as $\lambda\downarrow 0$
 exists and
        is uniform  in $\omega\in\Omega.$

 Assume that it is the limit of $V[\av_T]$ as $T\uparrow\infty.$
 It follows that \rref{slowlyT} and \rref{slowlyT_} hold for $U_T=V[v_T]$.
 From Propositions~\ref{av_bw} and~\ref{av_bw_}, we infer that, for all $\varepsilon>0$,
 $|V[\av_T](\omega)-V[\bw_{1/T}](\omega)|<\varepsilon$ holds for all $\omega\in\Omega$ if $T$ is  sufficient large.
 Thus,  the limit of $V[\bw_\lambda]$ as $\lambda\downarrow 0$ exists, is uniform, and coincides with the
 limit of $V[\av_T]$ as $T\uparrow\infty.$

  The case of the limit for $V[\bw_\lambda]$ is  considered analogously.
 It is only necessary to apply conditions
 \rref{slowlyl} and \rref{slowlyl_} and
  Propositions~\ref{bw_av} and \ref{bw_av_}.
\bo

{\bf Proof of Theorem~\ref{normal1}.}\
Let the set of all bounded maps  from $\mm{K}$ to $\mm{R}$
be the set $\mathfrak{C}$. This set
     satisfies \rref{conditions}, and
     $\bw_\lambda, \av_T$,
     $\zeta^{\avg_T}_{T,h},\xi^{\bwg_\lambda}_{\lambda,h}\in\mathfrak{C}$ for all positive $\lambda,T$ and natural $h$.
Note that the map $\cyg$ (see \rref{295}) takes each payoff $c:\mm{K}\to\mm{R}$  to a  function $\cyg[c]:\Omega\to\mm{R}$. Since
$\cyg$ satisfies conditions~\rref{conditions1}--\rref{conditions2}, this map
is a game value map.

To finish the proof,  we can now apply Theorem~\ref{normal2}.
\bo

\section{Tauberian theorem in discrete time setting}
\label{apply}

{\bf On dynamics.}\
  Let there be given sets $\Omega$, $\mm{K}$ and a running cost $g$, as before.
 Assume that we would like to consider any
  process  as a function  from $\{0,1,2,\dots\}$ to $\Omega$.
  Since, for such a function,  its definition can be completed in the form
 \begin{eqnarray}
  z(k+t)\rav z(k)\quad \forall k\in\{0,1,2,\dots\}, t\in(0,1), \label{579}
\end{eqnarray}
  we can propose that this function is from $\mm{R}_+$ to $\Omega$ and lies in $\mm{K}$.
  Thus, in this section, we can assume that all $z\in\mm{K}$ satisfy \rref{579}.

{\bf On payoffs.}\
Recall  that  $\mathfrak{U}$ is
   the set of all bounded maps from $\Omega$ to $\mm{R}$.
 Let us consider a non-empty set $\mathfrak{C}$ of  maps from $\mm{K}$ to $\mm{R}$
 and a game value map $V:\mathfrak{C}\to \mathfrak{U}$ satisfying conditions \rref{conditions} and
 \rref{conditions1},\rref{conditions2}, respectively.

 For all $\mu\in(0,1), n\in\mm{N}$, define payoffs $\bar{v}_n:\mm{K}\to\mm{R}, \bar{w}_\mu:\mm{K}\to\mm{R}$ as follows:
 \begin{eqnarray*}
 \bar{\av}_n(z)=\frac{1}{n}\sum_{t=0}^{n-1} g(z(t))\in [0,1],\quad
 \bar{\bw}_\mu(z)=\mu\sum_{t=0}^\infty (1-\mu)^{t} g(z(t))\in[0,1]\qquad \forall z\in{\mm{K}}.
\end{eqnarray*}
 Also, assume that
 $\bar{v}_n, \bar{w}_\mu\in\mathfrak{C}$  for all $\mu\in(0,1), n\in\mm{N}$.

 \begin{subequations}
{\bf On  Dynamic Programming Principle.}\
    Let us say that the family of payoffs  $\bar{\av}_n ({n\in\mm{N}})$
enjoys {\it the Dynamic Programming Principle}
iff, for all $n,h\in\mm{N}$,  the payoff
\begin{equation}\label{617z}
\bar{\mm{K}}\ni z\mapsto \frac{1}{n+h}\sum_{t=0}^{h-1} {g}(z(t))+\frac{n}{n+h}V[\bar{\av}_{n}](z(h))
\end{equation}
 lies in $\mathfrak{C}$ and the value of $V$  for this payoff
 coincides with $V[\bar{\av}_{n+h}]$

     Let us say that the family of payoffs  $\bar{\bw}_\mu (\mu>0)$
enjoys {\it the Dynamic Programming Principle} iff,  for all $\mu\in(0,1)$,  the function  $V[\bar{\bw}_\mu]$
 coincides with values of $V$  for payoffs
\begin{equation}\label{617x}
\mm{K}\ni z\mapsto
 \mu\sum_{t=0}^{h-1} (1-\mu)^{t} {g}(z(t))+(1-\mu)^{h}V[\bar{\bw}_\mu](z(h)) \quad\forall h\in\mm{N},
\end{equation}
and each of these payoffs lies in $\mathfrak{C}$.
 \end{subequations}

\begin{theorem}\label{normal3}
  Let all processes $z\in\mm{K}$ satisfy \rref{579}.
Let, for a game value map $V:\mathfrak{C}\to\mathfrak{U}$,
 the payoff families $\bar{\av}_n ({n\in\mm{N}})$ and  $\bar{\bw}_\mu (\mu>0)$
 enjoy the Dynamic Programming Principle.

Then, the following two statements are equivalent:
\begin{description}
  \item[$(\imath)$]   The sequence of functions $V[\bar{\av}_n]$ $(n\in\mm{N})$  converges uniformly on $\Omega$ as $n\uparrow\infty.$
  \item[$(\imath\imath)$] The family of functions $V[\bar{\bw}_\mu]$ $(0<\mu<1) $converges uniformly on $\Omega$ as $\mu\downarrow0$.
\end{description}
Moreover, when at least one of these statements holds, we have
    \begin{equation*}
    \lim_{n\uparrow\infty}V[\bar{\av}_n](\omega)=
       \lim_{\mu\downarrow 0}V[\bar{\bw}_\mu](\omega)\quad \forall\omega\in\Omega.
       \end{equation*}
\end{theorem}
{\bf Proof
 of Theorem~\ref{normal3}.}\
 We will use estimate~\rref{2012} and Lemma~\ref{765} proved in Appendix~\ref{AA}.

   Denote by $\mathfrak{C}_{b}$  the set of all bounded maps from $\mm{K}$ to $\mm{R}$.
   We can assume $\mathfrak{C}\subset\mathfrak{C}_{b}$; otherwise, we would always use
   $\mathfrak{C}\cap\mathfrak{C}_{b}$ instead of $\mathfrak{C}.$ Now,
   by  Lemma~\ref{765}, we can set  ${V}[c]\in\mathfrak{U}$ for all $c\in\mathfrak{C}_{b}\setminus\mathfrak{C}$
   such that  conditions \rref{conditions1}--\rref{conditions2} keep to hold.
   Thus, we obtain the game value map $V:\mathfrak{C}_{b}\to\mathfrak{U}.$

   Consider a function $\mu^*:\mm{R}_{>0}\to (0,1)$ defined as follows:
  $\mu^*(\lambda)=1-e^{-\lambda}$ for all $\lambda>0$.
  Then,
  $\int_{t}^{t+1}\lambda e^{-\lambda r}\,dr=e^{-\lambda t}\mu^*(\lambda)=
  (1-\mu^*(\lambda))^{t}\mu^*(\lambda)$
  for all $t\geq 0,\lambda>0$.
  Now, \rref{579} implies $\bar{\bw}_{\mu^*(\lambda)}\equiv \bw_{\lambda}$
  for all $\lambda>0$.   Note that $\mu^*(0+)=0+.$
  Then,  the limit of $V[\bar{\bw}_\mu]$ as $\mu\downarrow 0$ exists and is uniform in $\Omega$ iff
   the limit of $V[{\bw}_\lambda]$ as $\lambda\downarrow 0$ exists and is uniform in $\Omega$.
   Moreover, in this case, these limits coincide.

  Also, it is easy to see  that $\bar{\av}_n\equiv \av_n$ for all $n\in\mm{N}.$
 Since the payoffs $\av_{T}(T>0)$ are bounded, we have $\av_{T}\in\mathfrak{C}_b$ for all $T>0.$
    For each $T>0$, we can choose $n\in\mm{N}$ such that $T\in(n-1,n]$, hence, we obtain
$$
  |V[\bar{\av}_{n}](\omega)-V[{\av}_{T}](\omega)|\leqref{2012} \frac{2(n-T)}{T}\leq\frac{2}{T}
  \qquad\forall \omega\in\Omega.$$
  Thus,  the limit of $V[\bar{\av}_T]$ as  $T\uparrow \infty$ exists and is uniform in $\Omega$ iff
   the limit of $V[\bar{\av}_n]$ as $n\to\infty$ exists and is uniform in $\Omega$.
   Moreover, in this case, these limits coincide.

  At last, thanks to \rref{579} and $\bar{\bw}_{\mu^*(\lambda)}\equiv \bw_{\lambda}$, $\bar{\av}_n\equiv \av_n$ for all $n\in\mm{N},\lambda>0,$ we have that, for all natural $T=n,h\in\mm{N}$ and positive $\lambda,$
   payoffs \rref{617z} and \rref{617x} coincide with $\zeta^{V[\av_T]}_{T,h}$
and  $\xi^{V[\bw_\lambda]}_{\lambda,h}$, respectively.
     Therefore, the payoff families  $\bw_\lambda (\lambda>0)$ and $\av_T (T>0)$
enjoy the weak Dynamic Programming Principle.

 We have verified all conditions of Theorem~\ref{normal2}. Moreover, the corresponding limits in Theorem~\ref{normal2} and in Theorem~\ref{normal3} exist and are uniform   in
 $\omega\in\Omega$ only simultaneously. Thanks to Theorem~\ref{normal2}, Theorem~\ref{normal3} is proved.
\bo

Let us showcase the application of this approach in a stochastic framework.

  Consider a $\sigma$-algebra $\ct{A}$ on $\mm{K}$.
     Assume that $\mathfrak{C}_{\ct{A}}$ is the set of all $\ct{A}$-measurable bounded maps of $\mm{K}$ to $\mm{R}$. It is easy to see that $\mathfrak{C}_\ct{A}$ satisfies the condition \rref{conditions}.

  Similarly to Section~\ref{abstractgame}, for all $\omega\in\Omega$, let there be given non-empty sets $\ct{L}(\omega)$ and $\ct{M}(\omega)$.
  Let, for all   $\omega\in\Omega$, each pair $(l,m)\in\ct{L}(\omega)\times\ct{M}(\omega)$ of players' rules
  induce  a probability distribution $\mm{P}^\omega_{lm}$ (over $(\mm{K},\ct{A})$) along with its
   mathematical expectation $\mm{E}^\omega_{lm}$.
  Similarly to \rref{295}, define the map $W:\mathfrak{C}_{\ct{A}}\to \mathfrak{U}$ by the following rule: for each $c\in\mathfrak{C}$,
  \begin{eqnarray}
 W[c](\omega)\rav\sup_{l\in \ct{L}(\omega)}\inf_{m\in \ct{M}(\omega)} \mm{E}^\omega_{lm} c\qquad \forall\omega\in\Omega.
 \label{295_}
\end{eqnarray}
Evidently, $W:\mathfrak{C}\to \mathfrak{U}$ satisfies conditions \rref{conditions1},\rref{conditions2}. Thus, $W$ is a game value map.

 Applying Theorem~\ref{normal3} for this game value map $W,$ we obtain
\begin{corollary}\label{normal11}
  Let all processes $z\in\mm{K}$ satisfy \rref{579}.  Also, assume that
   for all $\mu\in(0,1)$, $n\in\mm{N}$ the payoffs
 $\bar{v}_n, \bar{w}_\mu$ are $\ct{A}$-measurable.

Let, for the game value map $W:\mathfrak{C}\to\mathfrak{U}$ (see \rref{295_}),
 the payoff families $\bar{\av}_n ({n\in\mm{N}})$ and  $\bar{\bw}_\mu (0<\mu<1)$
 enjoy the Dynamic Programming Principle.

 Then, the following two statements are equivalent:
\begin{description}
  \item[$(\imath)$]   The maps $\displaystyle\Omega\ni\omega\mapsto\sup_{l\in \ct{L}(\omega)}\inf_{m\in \ct{M}(\omega)} \mm{E}^\omega_{lm} \bar{\av}_n$ 
      converge uniformly on $\Omega$ as $n\uparrow\infty.$
  \item[$(\imath\imath)$] The maps $\displaystyle\Omega\ni\omega\mapsto\sup_{l\in \ct{L}(\omega)}\inf_{m\in \ct{M}(\omega)} \mm{E}^\omega_{lm} \bar{\bw}_\mu$ 
       converge uniformly on $\Omega$ as $\mu\downarrow0$.
\end{description}
Moreover, when at least one of these statements holds, we have
    \begin{equation*}
    \lim_{n\uparrow\infty}
    \sup_{l\in \ct{L}(\omega)}\inf_{m\in \ct{M}(\omega)} \mm{E}^\omega_{lm} \bar{\av}_n=
    \lim_{\mu\downarrow 0}
        \sup_{l\in \ct{L}(\omega)}\inf_{m\in \ct{M}(\omega)} \mm{E}^\omega_{lm} \bar{\bw}_\mu\quad \forall\omega\in\Omega.
       \end{equation*}
\end{corollary}

\section{Differential games without saddle point.}
\label{diff}

{\bf Dynamic equation.}\
 Consider a nonlinear system in $\mm{R}^m$ controlled by two players,
\begin{equation}\label{sys}
 \dot{x}=f(x,a,b),\ x(0)\in \mm{R}^m, a(t)\in \mm{A},\ b(t)\in \mm{B};
\end{equation}
here,~$\mm{A}$ and~$\mm{B}$ are non-empty compact subsets of finite-dimensional Euclidean spaces.

 In this section, we assume that
\begin{enumerate}
  \item the functions
 $f:\mm{R}^m \times \mm{A} \times \mm{B}\to\mm{R}^m$, $g:\mm{R}^m \times \mm{A} \times \mm{B}\to [0,1]$ are continuous;
  \item these functions are Lipschitz continuous in the state variable, namely, for a constant $L>0$,
  $$\big|\big|f(x,a,b)-f(y,a,b)\big|\big|+\big|g(x,a,b)-g(y,a,b)\big|\leq L\big|\big|x-y\big|\big|
  \ \forall x,y\in\mm{R}^m,a\in \mm{A},b\in \mm{B}.$$
\end{enumerate}

Denote by  $B(\mm{R}_+,\mm{A})$ and by $B(\mm{R}_+,\mm{B})$ the sets of all Borel measurable functions $\mm{R}_+\ni t\mapsto a(t)\in
\mm{A}$ and $\mm{R}_+\ni t\mapsto b(t)\in \mm{B}$, respectively.
 Now, for each pair  $(a,b)\in B(\mm{R}_+,\mm{A})\times B(\mm{R}_+,\mm{B})$, for every initial condition  $x(0)=x_*$, system
 \rref{sys} generates the unique solution  $x(\cdot)=y(\cdot;x_*,a,b)$ defined
 for the whole
 $\mm{R}_+$. Denote by $Y(x_*)$ the set of all such solutions with $x(0)=x_*$.

  Consider a set $\mm{X}\subset\mm{R}^m$ that is strongly invariant with respect to system \rref{sys}, i.e.,
   let $x(t)\in\mm{X}$ for all $t\in\mm{R}_+$, $x_*\in\mm{X},x\in Y(x_*)$. Set
  $\mm{Y}\rav \cup_{x_*\in\mm{X}} Y(x_*).$

{\bf On strategies of players.}\
    Like before, let the goal of the first player be to maximize the payoff function and let the task of the second one be to
    minimize it.
      Our payoff functions are as follows: for all $\lambda,T>0,(x,a,b)\in\mm{Y}\times B(\mm{R}_+,\mm{A})\times B(\mm{R}_+,\mm{B}),$
     \begin{eqnarray*}
 \av_T(x,a,b)\rav\frac{1}{T}\int_{0}^T g(x(t),a(t),b(t))\,dt,\\
 \bw_\lambda(x,a,b)\rav\lambda\int_{0}^\infty e^{-\lambda t} g(x(t),a(t),b(t))\,dt.
 \end{eqnarray*}

 In the general case, without Isaacs's condition,
  the lower game value and the upper game value  depend on  choosing the formalization of strategies of players
  \cite[Ch.~XVI]{ks}, \cite[Subsect.~14]{subb}, \cite{Bardi}.
     For simplicity, we now assume that the first player announces a  nonanticipating  strategy
      (see \cite{EK,Rox,RN}) and another, knowing it, selects an admissible measurable control.
\begin{definition}
      A map  $\alpha : B(\mm{R}_+,\mm{B})\mapsto B(\mm{R}_+,\mm{A})$ is called  {\it a  nonanticipating strategy}   for  the  first
      player   if,  for  all  $t>0$ and $b,b'\in B(\mm{R}_+,\mm{B})$,
      $b|_{[0,t]}   =  b'|_{[0,t]}$   implies that $\alpha[b]|_{[0,t]}  =  \alpha[b']|_{[0,t]}$.
  We  denote  by $\ct{A}$  the  set  of all nonanticipating  strategies  for  the  first  player.
\end{definition}
   Now, for all $T,\lambda>0,$ we can define lower game values,
     \begin{eqnarray*}
 \avg_T(x_*)&\rav&\sup_{\alpha\in\ct{A}}\inf_{b\in B(\mm{R}_+,\mm{B})}
   \av_T(y(\cdot;x_*,\alpha(b),b),\alpha(b),b)\qquad\forall x_*\in\mm{X},\\
 \bwg_\lambda(x_*)&\rav&\sup_{\alpha\in\ct{A}}\inf_{b\in B(\mm{R}_+,\mm{B})}
   \bw_\lambda(y(\cdot;x_*,\alpha(b),b),\alpha(b),b)\qquad\forall x_*\in\mm{X}.
 \end{eqnarray*}

{\bf Constructing a game value map.} \
 Let us set
    $$\Omega\rav\mm{X}\times\mm{A}\times\mm{B},\quad \mm{K}\rav \mm{Y}\times B(\mm{R}_+,\mm{A})\times B(\mm{R}_+,\mm{B}).$$
     Let ${\mathfrak{C}}$ be the set of all bounded  maps from $\mm{K}$ to $\mm{R}$; this set
     satisfies \rref{conditions}, and
     $\bw_\lambda, \av_T$,
     $\zeta^{\avg_T}_{T,h},\xi^{\bwg_\lambda}_{\lambda,h}\in\mathfrak{C}$ for each positive $\lambda,T$ and natural $h$.
For  every map $c\in{\mathfrak{C}}$,  we can consider the following value function:
     \begin{eqnarray*}
   V[c](x_*,a_*,b_*)\rav\sup_{\alpha\in\ct{A}}\inf_{b\in B(\mm{R}_+,\mm{B})}
   c(y(\cdot;x_*,\alpha(b),b),\alpha(b),b)\qquad\forall \omega=(x_*,a_*,b_*)\in\Omega.
\end{eqnarray*}
      It is easy prove that $V:\mathfrak{C}\to\mathfrak{U}$ satisfies \rref{conditions1},\rref{conditions2}. Thus, we obtain the
      game value map $V.$

The dynamic programming principle with respect to  nonanticipating strategies  for  Bolza
functionals (particularly, for the payoffs $v_T$) is well-known, see \cite{evans,subb}. Such a principle for the
payoff function $\bw_\lambda$ follows from \cite[Theorem VIII.1.9]{BardiGame}. All conditions of Theorem~\ref{normal2} verify.
Moreover, $V[\bw_\lambda]\equiv\bwg_\lambda$, $V[\av_T]\equiv\avg_T$ are independent of $a_*,b_*$. Now,
   thanks to Theorem~\ref{normal2}, we obtain
\begin{theorem}
     \label{maintheoremdiff}
  Let $f,g$ be as before, and let a non-empty set $\mm{X}\subset\mm{R}^m$  be strongly invariant with respect to  \rref{sys}.

  Then,    for a function $U_*:\mm{X}\to [0,1]$,
 $U_*$ is  a limit of $\avg_T$ as $T\uparrow \infty$
 that is uniform in $\mm{X}$ iff
  $U_*$ is  a  limit of $\bwg_\lambda$ as $\lambda\downarrow 0$
 that is uniform in $\mm{X}.$
    \end{theorem}
   Under Isaacs's condition, this theorem was proved in \cite{Khlopin3}.

  \section*{Acknowledgements}
   I would like to
express my gratitude to  Ya.V. Salii for
the translation. 
I am grateful to an anonymous referee for helpful comments.
This study was supported by the Russian Foundation for
the Program of Basic Researches of the Ural Branch of the Russian Academy of Sciences  (Project 15-16-1-8 ``Control Synthesis under Incomplete Data, the Reachability Problems and Dynamic Optimization'') and by  by the Russian Foundation for Basic Research, project no.
16-01-00505.

\appendix

\section{Auxiliary statements.}
\label{AA}
\begin{subequations}
  Assume a non-empty set $\mathfrak{C}$ satisfies \rref{conditions} and a game value map $V:\mathfrak{C}\to\mathfrak{U}$ enjoys
    \rref{conditions1},  \rref{conditions2}.

{\bf Two estimates.}\
  Consider  $T>0,r>1$, and $z\in\mm{K}.$ Assume $\av_T,\av_{rT}\in\mathfrak{C}$. Now,
\begin{eqnarray*}
  |\av_T(z)-\av_{rT}(z)|&=&
  \Big|\frac{1}{T}\int_0^{T} g(z(t))\,dt-\frac{1}{rT}\int_0^{rT} g(z(t))\,dt\Big|\\
  &\leq& \Big|\frac{1}{T}-\frac{1}{rT}\Big|\int_0^T g(z(t))\,dt+\frac{1}{rT} \int_T^{rT}g(z(t))\,dt\\
  &\leq&
  1-\frac{1}{r}+\frac{r-1}{r}=\frac{2r-2}{r}<2(r-1).
\end{eqnarray*}
 So, $\av_{rT}(z)-2(r-1)\leq \av_T(z)\leq\av_{rT}(z)+2(r-1)$ for all $z\in\mm{K}$.  Thanks to
 \rref{conditions1},\rref{conditions2}, we obtain
\begin{eqnarray}\label{2012}
\Big|V[\av_T](\omega)-V[\av_{rT}](\omega)\Big|\leq 2(r-1)\qquad \forall\omega\in\Omega.
\end{eqnarray}

  Consider  $\lambda>0,r>1$, and $z\in\mm{K}.$  Assume $\bw_\lambda,\bw_{r\lambda}\in\mathfrak{C}$. Now,
\begin{eqnarray*}
  |\bw_\lambda(z)-\bw_{r\lambda}(z)|&=&
  \Big|\lambda\int_0^{\infty} e^{-\lambda t} g(z(t))\,dt-
  r\lambda\int_0^{\infty}e^{-r\lambda t} g(z(t))\,dt\Big|\\
  &\leq&
  \int_0^{\infty} |\lambda e^{-\lambda t}-r\lambda e^{-r\lambda t}|g(z(t))\,dt\\
  &\leq&
  \int_0^{\infty} |\lambda e^{-\lambda t}-r\lambda e^{-r\lambda t}|\,dt\\
   &=&
  2\int_0^{\infty} \max\{\lambda e^{-\lambda t},r\lambda e^{-r\lambda t}\}\,dt-\int_0^{\infty} \lambda e^{-\lambda
  t}\,dt-\int_0^{\infty} r\lambda e^{-r\lambda t}\,dt\\
  &\leq&
  2\int_0^{\infty} r\lambda e^{-\lambda t}\,dt-2=2(r-1).
\end{eqnarray*}
 Thus,
\begin{eqnarray}\label{2022}
\Big|V[\bw_\lambda](\omega)-V[\bw_{r\lambda}](\omega)\Big|\leq 2(r-1)
\qquad \forall\omega\in\Omega.
\end{eqnarray}
\end{subequations}

{\bf On extension of $V$.}\
   Denote by $\mathfrak{C}_{b}$   the set of all bounded  maps from $\mm{K}$ to $\mm{R}$.
   Assume $\mathfrak{C}\subset\mathfrak{C}_b$.

Note that \rref{conditions} holds for $\mathfrak{C}_{b}$.
Define a map $\mm{V}_b:\mathfrak{C}_{b}\to\mathfrak{U}$ as follows:
$$
 \mm{V}_b[c'](\omega)\rav \sup\{V[c](\omega)\,|\,c\in\mathfrak{C}, c\leq c'\}
 \qquad \forall \omega\in{\Omega},c'\in\mathfrak{C}_b.
  $$
  Since, for all $R\in\mm{R},c\in\mathfrak{C}$, the inequality $-R\leq c\leq R$ implies $ -R=V[-R]\leq V[c]\leq V[R]=R$ by \rref{conditions1},\rref{conditions2}, the function $\mm{V}_b[c']$ is well-defined and  bounded for every $c'\in\mathfrak{C}_{b}$.
 It is easy to prove that \rref{conditions1} and  \rref{conditions2} hold for $\mm{V}_b.$
 In addition,
$$
 \mm{V}_b[c''](\omega)\rav \sup\{V[c](\omega)\,|\,c\in\mathfrak{C}, c\leq c''\}\ravref{conditions2}V[c''](\omega)
 \qquad \forall \omega\in{\Omega},c''\in\mathfrak{C}.
  $$
  So, $\mm{V}_b|_{\mathfrak{C}}\equiv V.$
 Thus, we have proved
\begin{lemma}\label{765}
   Let a non-empty set $\mathfrak{C}\subset\mathfrak{C}_b$ satisfy \rref{conditions} and a map $V:\mathfrak{C}\to\mathfrak{U}$
   enjoy
    \rref{conditions1} and  \rref{conditions2}.
  Then,  the map $\mm{V}_b:\mathfrak{C}_{b}\to\mathfrak{U}$ enjoys \rref{conditions1} and  \rref{conditions2} and
  is an extension
   of $V$.
\end{lemma}

\section{ The proof of Proposition \ref{av_bw}.}
\label{A}

{\bf On simplicity of notation.}\

For all payoff functions $c$, let us also use the following notation:
$$\ssp{c}{\omega}\equiv V[c](\omega)\qquad\forall\omega\in\Omega.$$
For instance, for  a function $U_*:\Omega\to\mm{R}$ and $\omega\in\Omega$,  the symbol $\ssp{U_*(z(1))}{\omega}$ means $V[c](\omega)$ for the payoff
function  $\mm{K}\ni z\mapsto c(z)=U_*(z(1))\in\mm{R}$. The symbol  $\ssp{U_*(z_1(1))}{\omega}$ also means the same
for the same payoff.
Moreover, the symbols
$$ \sspk{\int_{0}^h a(t)g(z(t))\,dt+ U_*(z(h))}{\omega},\
\ssp{\ssp{U_*(z_1(1))}{z(1)}}{\omega},\
 \sspk{\int_{0}^h b(t)g(z_1(t))\,dt+ \ssp{c}{z_1(h)}}{z(h')}$$
 denote the value of the game value map   for the payoff function
 $\mm{K}\ni z\mapsto \int_{0}^h a(t)g(z(t))\,dt+ U_*(z(h))\in\mm{R}$ at $\omega$,
 for the payoff function
 $\mm{K}\ni z\mapsto c_1(z)=\ssp{U_*(z_1(1))}{z(1)}\in\mm{R}$  at $\omega$,
 for  the payoff function
 $$\mm{K}\ni z_1\mapsto \int_{0}^h b(t)g(z_1(t))\,dt+V[c](z_1(h))\in\mm{R}$$ at $z(h')$,
 respectively.
 Now, the expression $\ssp{\ssp{U_*(z_5(1))}{z_{2}(1)}}{\omega}$
  is equal to $\ssp{\ssp{U_*(z_1(1))}{z(1)}}{\omega}$,  similarly to e.g. the equivalence of $\int_A f(x)\int_B g(y)dydx$ and $\int_Af(r)\int_B g(s)dsdr.$

{\bf Auxiliary estimates.}\
 Recall that
 $\mathfrak{C}_{b}$ is the set of all bounded maps from $\mm{K}$ to $\mm{R}$.
 Without loss of generality, we  assume that $\mathfrak{C}\subset\mathfrak{C}_{b}$; otherwise, we could always use
   $\mathfrak{C}\cap\mathfrak{C}_{b}$ instead of $\mathfrak{C}.$ Further,
   we can assume that  $\mathfrak{C}=\mathfrak{C}_{b}$;
   otherwise,  applying Lemma~\ref{765}, we could always define  $V$ on
   $\mathfrak{C}_{b}\setminus \mathfrak{C}$
   such that  conditions \rref{conditions1}--\rref{conditions2} keep to hold.
   Thus, we obtain the value $V[c]$ for all bounded payoffs $c:\mm{K}\to\mm{R}$.

Set
$$\kappa(T,p_0)\rav \sup_{p\in[1,p_0]}\sup_{\omega\in\Omega}\,
\big(\,U_{Tp^{-1}}(\omega)-U_{T}(\omega)\big)\quad\forall T>0,p_0>1.$$

Fix some positive $\epsi<1/4$. Applying $2\epsi\leq \epsi^{1/2}$,
we can choose some natural  $k\geq 2$ such that  $\epsi^{1/2}\leq k\epsi\leq (k+1)\epsi\leq 2{\epsi}^{1/2}$.
Set $p\rav 1+\frac{1}{k}$; we have
$$\frac{\epsi}{p-1}\leq\frac{p\epsi}{p-1}=(k+1)\epsi\leq 2{\epsi}^{1/2},\quad  2(p-1)=\frac{2}{k}\leq 2{\epsi}^{1/2}.$$

  By the definitions of the subsolution and the supersolution, for each $\epsi>0$,
    there exists a positive  $\bar{T}$ such that, for all $\lambda\in(0,1/\bar{T}),\omega\in\Omega$,
    and natural $h>\bar{T}$,$T>\bar{T}+h$,
\begin{eqnarray}
 U_{T}(\omega)&\leq&V\big[\zeta^{U_{T-h}}_{h,T}\big](\omega)+\epsi/2=
 \ssp{\frac{1}{T}\int_{0}^{h}g(z(t))\,dt+\frac{T-h}{T} U_{T-h}(z(h))+\epsi/2}{\omega},\label{14421}\\
 V[\bw_{\lambda}](\omega)&\geq&V\big[\xi^{V[\bw_{\lambda}]}_{h,\lambda}\Big](\omega)-\epsi/3=\ssp{\int_{0}^{h}\lambda e^{-\lambda t}g(z(t))\,dt+e^{-\lambda h}
 V[\bw_{\lambda}](z(h))-\epsi/3}{\omega}. \label{1442}
\end{eqnarray}
By \rref{slowlyT},  there exists  $\hat{T}>\max\{\bar{T},\frac{1}{\ln p}\}$ such that
 $\kappa(T,p')\leq\kappa(T,p)\leq\epsi/2$ holds
 for all $T>\hat{T}$, $p'\in(1,p]$.

{\bf The special choice of $\lambda$.}\
Fix every $\lambda<\frac{1}{(k+1)\hat{T}}$ such that $\frac{\ln p}{\lambda}$ is a natural number; set
\begin{eqnarray}\label{1444}
 h\rav\frac{\ln p}{\lambda},\ T\rav\frac{p\ln p}{\lambda(p-1)}=(k+1)h,\ q\rav p^{-1}=e^{-\lambda h}
 =\frac{T-h}{T},\ \end{eqnarray}
Observe that $h$, $T=(k+1)h,$ and $T-h=kh$ are natural.
It is easy to verify that $p>1$ implies that $\ln p<p-1< p\ln p$.
Then,
\begin{eqnarray}\label{1445}
 \frac{p}{\lambda}>T=\frac{p\ln p}{\lambda(p-1)}=\frac{h}{1-q}>\frac{1}{\lambda}>(k+1)\hat{T}.
\end{eqnarray}
 In particular,  from the inequalities $T=h(k+1)>(k+1)\hat{T}$ and $T-h\geq h>\hat{T}$, it follows
  that $T,h,\lambda$ satisfy \rref{14421},\rref{1442}. Also,
 $T\lambda=\frac{p\ln p}{p-1}\in(1,p)$ guarantees that
\begin{eqnarray}\label{1446}
 U_{1/\lambda}(\omega)\leq U_{T}(\omega)+\kappa(T,T\lambda)\leq  U_{T}(\omega)+\epsi/2
\quad\forall \omega\in\Omega.
\end{eqnarray}

{\bf Forward-tracking.}\
For all $\omega\in\Omega,$ from \rref{14421} it follows that
\begin{eqnarray}
 {U_{T}(\omega)} &\leqref{14421}&
 \ssp{\frac{1}{T}\int_{0}^{h}g(z(t))\,dt+\frac{T-h}{T} U_{T-h}(z(h))+\epsi/2}{\omega}\nonumber\\
 &\leq&\ssp{\frac{1}{T}\int_{0}^{h}g(z(t))\,dt+\frac{T-h}{T} U_{T}(z(h))+\epsi/2+\kappa(T,p)}{\omega}\nonumber\\
 &\leq&\ssp{\frac{1}{T}\int_{0}^{h}g(z(t))\,dt+q U_{T}(z(h))}{\omega}+\epsi.\label{144210}
\end{eqnarray}
The inequality will still hold if we replace the symbol~$z$ inside the square brackets with e.g.~$z_1$.
 We have
\begin{eqnarray*}
 {U_{T}(\omega)} &\leqref{144210}&\ssp{\frac{1}{T}\int_{0}^{h}g(z_1(t))\,dt+q U_{T}(z_1(h))}{\omega}+\epsi.
\end{eqnarray*}
In particular, for all $z\in\mm{K}$, we get
\begin{eqnarray*}
 {U_{T}(z(h))} &\leqref{144210}&\ssp{\frac{1}{T}\int_{0}^{h}g(z_1(t))\,dt+q U_{T}(z_1(h))}{z(h)}+\epsi.
\end{eqnarray*}
 Substituting the corresponding part into \rref{144210}, we obtain
 \begin{eqnarray*}
 {U_{T}(\omega)}
 &\leq&\ssp{\frac{1}{T}\int_{0}^{h}g(z(t))\,dt+q U_{T}(z(h))}{\omega}+\epsi\\
 &\leqref{conditions2}&\ssp{\frac{1}{T}\int_{0}^{h}g(z(t))\,dt+
 \ssp{\frac{1}{T}\int_{0}^{h}q g(z_1(t))\,dt+q^{2} U_{T}(z_1(h))}{z(h)}+q\epsi}{\omega}+\epsi\\
 &\ravref{conditions1}&\ssp{\frac{1}{T}\int_{0}^{h}g(z(t))\,dt+
 \ssp{\frac{1}{T}\int_{0}^{h}q g(z_1(t))\,dt+q^{2} U_{T}(z_1(h))}{z(h)}}{\omega}+\epsi(1+q).\\
\end{eqnarray*}
 Repeating, we obtain
 \begin{eqnarray*}
 {U_{T}(\omega)}
 &\leq&\sspll{\frac{1}{T}\int_{0}^{h}g(z(t))\,dt+
 \sspl{\frac{1}{T}\int_{0}^{h}q g(z_1(t))\,dt}\\
 &\ &\sspr{+\sspk{\frac{1}{T}\int_{0}^{h}q^2 g(z_2(t))\,dt+q^{3}
 U_{T}(z_2(h))}{z_1(h)}+q^2\epsi}{z(h)}}\ssprr{}{\omega}+\epsi(1+q)\\
 &=&\sspll{\frac{1}{T}\int_{0}^{h}g(z(t))\,dt+
 \sspl{\frac{1}{T}\int_{0}^{h}q g(z_1(t))\,dt}\\
 &\ &\sspr{+\sspk{\frac{1}{T}\int_{0}^{h}q^2 g(z_2(t))\,dt+q^{3} U_{T}(z_2(h))}{z_1(h)}}{z(h)}}\ssprr{}{\omega}+\epsi(1+q+q^2).
\end{eqnarray*}
Proceeding in a similar way,   for all $n\in\mm{N}$, $\omega\in\Omega$,  we obtain
 \begin{eqnarray*}
 {U_{T}(\omega)}
 &\leq&\ssplll{\frac{1}{T}\int_{0}^{h}g(z(t))\,dt+
 \sspll{\frac{1}{T}\int_{0}^{h}q g(z_1(t))\,dt+ \sspl{\frac{1}{T}\int_{0}^{h}q^2 g(z_2(t))\,dt}}}\\
 &\ &\ssprrr{\ssprr{\sspr{+\dots+\sspk{\frac{1}{T}\int_{0}^{h}q^n g(z_n(t))\,dt+q^{n+1}
 U_{T}(z_n(h))}{z_{n-1}(h)}\dots}{z_1(h)}}{z(h)}}{\omega}\\
 &\ &+\epsi(1+q+q^2+\dots+q^n).
\end{eqnarray*}

{\bf Backtracking.}\
 Since $U_{T}$ is bounded, we can
 choose natural $n$ such that $U_T q^{n+1}\leq \epsi^{1/2}$.
 Then,
 \begin{eqnarray*}
 {U_{T}(\omega)}
 &\leq&\ssplll{\frac{1}{T}\int_{0}^{h}g(z(t))\,dt+
 \sspll{\frac{1}{T}\int_{0}^{h}q g(z_1(t))\,dt+ \sspl{\frac{1}{T}\int_{0}^{h}q^2 g(z_2(t))\,dt}}}\\
 &\ &\ssprrr{\ssprr{\sspr{+\dots+\sspk{\frac{1}{T}\int_{0}^{h}q^n
 g(z_n(t))\,dt+\epsi^{1/2}}{z_{n-1}(h)}\dots}{z_1(h)}}{z(h)}}{\omega}+\frac{\epsi}{1-q}\\
 &\leq&\ssplll{\frac{1}{T}\int_{0}^{h}g(z(t))\,dt+
 \sspll{\frac{1}{T}\int_{0}^{h}q g(z_1(t))\,dt+ \sspl{\frac{1}{T}\int_{0}^{h}q^2 g(z_2(t))\,dt}}}\\
 &\ &\ssprrr{\ssprr{\sspr{+\dots+\sspk{\frac{1}{T}\int_{0}^{h}q^n
 g(z_n(t))\,dt}{z_{n-1}(h)}\dots}{z_1(h)}}{z(h)}}{\omega}+\epsi^{1/2}+\frac{\epsi}{1-q}.
\end{eqnarray*}
 By the choice of $p$, we have $\frac{\epsi}{1-q}=\frac{p\epsi}{p-1}\leq 2{\epsi}^{1/2}.$
 In addition, $T\lambda\geq 1$ by \rref{1445}. Now, $g\geq0$ leads to
\begin{eqnarray*}
U_{T}(\omega) &\leq&\ssplll{\frac{1}{T}\int_{0}^{h}g(z(t))\,dt+
 \sspll{\frac{1}{T}\int_{0}^{h}q g(z_1(t))\,dt+ \sspl{\frac{1}{T}\int_{0}^{h}q^2 g(z_2(t))\,dt}}}\\
 &\ &\ssprrr{\ssprr{\sspr{+\dots+\sspk{\frac{1}{T}\int_{0}^{h}q^n
 g(z_n(t))\,dt}{z_{n-1}(h)}\dots}{z_1(h)}}{z(h)}}{\omega}+\epsi^{1/2}+2{\epsi}^{1/2}\\
 &\leq&\ssplll{\int_{0}^{h}\lambda g(z(t))\,dt+
 q\sspll{\int_{0}^{h}\lambda g(z_1(t))\,dt+ q\sspl{\int_{0}^{h} \lambda g(z_2(t))\,dt}}}\\
 &\ &\ssprrr{\ssprr{\sspr{+\dots+q\sspk{\int_{0}^{h}\lambda
 g(z_n(t))\,dt}{z_{n-1}(h)}\dots}{z_1(h)}}{z(h)}}{\omega}+3{\epsi}^{1/2}.
\end{eqnarray*}
   Recall that  $\lambda=\frac{\ln p}{h}$ by \rref{1444}. Also, thanks to \rref{1444}, we have
   $pe^{-\lambda t}\geq pe^{-\lambda h}=pq=1$  for all $t\in[0,h]$. It follows from $g\geq 0$ and $V[\bw_\lambda]\geq 0$ that
\begin{eqnarray*}
U_{T}(\omega) &\leq&\ssplll{\int_{0}^{h}p\lambda e^{-\lambda t}g(z(t))\,dt+
 q\sspll{\int_{0}^{h}p\lambda e^{-\lambda t}g(z_1(t))\,dt}}\\
 &\ &\ssprrr{\ssprr{+\dots+q\sspk{\int_{0}^{h}p\lambda e^{-\lambda t}
 g(z_n(t))\,dt+V[\bw_\lambda](z_{n}(h))}{z_{n-1}(h)}\dots}{z(h)}}{\omega}+3{\epsi}^{1/2}\\
 &\ravref{conditions1}&p\ssplll{\int_{0}^{h}\lambda e^{-\lambda t}g(z(t))\,dt+
 q\sspll{\int_{0}^{h}\lambda e^{-\lambda t}g(z_1(t))\,dt}}\\
 &\ &\ssprrr{\ssprr{+\dots+q\sspk{\int_{0}^{h}\lambda e^{-\lambda t}
 g(z_n(t))\,dt+qV[\bw_\lambda](z_{n}(h))}{z_{n-1}(h)}\dots}{z(h)}}{\omega}+3{\epsi}^{1/2}.
\end{eqnarray*}
 Since  $V[\bw_\lambda]$ is a supersolution   (see \rref{1442}), in view of $q=e^{-\lambda h}$, we obtain
 $$q\sspk{\int_{0}^{h}\lambda e^{-\lambda t} g(z_n(t))\,dt+e^{-\lambda h}V[\bw_\lambda](z_{n}(h))}{z_{n-1}(h)}\leq e^{-\lambda
 h}V[\bw_\lambda](z_{n-1}(h))+q\epsi/3.$$
Thus,
\begin{eqnarray*}
U_{T}(\omega)
 &\leq&p\ssplll{\int_{0}^{h}\lambda e^{-\lambda t}g(z(t))\,dt+
 q\sspll{\int_{0}^{h}\lambda e^{-\lambda t}g(z_1(t))\,dt}}
 +\dots+q\sspl{\int_{0}^{h}\lambda e^{-\lambda t} g(z_{n-1}(t))\,dt}\\&\ &+\ssprrr{\ssprr{\sspr{e^{-\lambda
 h}V[\bw_\lambda](z_{n-1}(h))}{z_{n-2}(h)}\dots}{z(h)}}{\omega}+pq^{n}\epsi/3+3{\epsi}^{1/2}.
\end{eqnarray*}
Proceeding in a similar way,  we have
\begin{eqnarray*}
U_{T}(\omega) &\leq&
 p\ssp{\int_{0}^{h}\lambda e^{-\lambda t} g(z(t))\,dt+e^{-\lambda h}V[\bw_\lambda](z(h))}{\omega}
 +p(q+q^2+\dots+q^{n})\epsi/2+3{\epsi}^{1/2}\\
 &\leq&p V[\bw_\lambda](\omega)+p(1+q+q^2+\dots+q^{n})\epsi/3+3{\epsi}^{1/2}\\
&\leq&
 p V[\bw_\lambda](\omega)+2p{\epsi}^{1/2}/3+3{\epsi}^{1/2}.
\end{eqnarray*}
Using $p\leq 1+\frac{1}{2}$ and $(p-1)V[\bw_\lambda]\leq p-1\leq {\epsi}^{1/2}$, we obtain
\begin{eqnarray}\nonumber
U_{1/\lambda}\leqref{1446}U_{T}+\epsi/2&\leq&
 p V[\bw_\lambda]+5{\epsi}^{1/2}+\epsi/2\\
 &\leq& V[\bw_\lambda]+(p-1)+5{\epsi}^{1/2}+\epsi/2\leq
 V[\bw_\lambda]+6{\epsi}^{1/2}+\epsi/2\label{2250}
\end{eqnarray}
 for all positive $\lambda<\frac{1}{(k+1)\hat{T}}$ such that $\frac{\ln p}{\lambda}$ is a natural number.

{\bf The general case.}\
 For positive $\lambda'<\frac{1}{(k+1)\hat{T}}$, we can choose positive $r>1$ such that $\frac{r\ln p}{\lambda'}$  is a natural
 number
  and $0\leq \frac{r\ln p}{\lambda'}-\frac{\ln p}{\lambda'}\leq 1$.  Recall that, by the choice of $\hat{T}$, we have  $\hat{T}
  \ln p\geq 1.$
  By the choice of $\hat{T}$ and $\lambda',$ we obtain
 $$1\leq r\leq 1+\frac{\lambda'}{\ln p}\leq  1+\frac{1}{k+1}<p.$$
  First,  thanks to  \rref{2022},  we get
   $$V[\bw_{\lambda'/r}]\leqref{2022} V[\bw_{\lambda'}]+2(r-1)\leq V[\bw_{\lambda'}]+2(p-1)\leq
   V[\bw_{\lambda'}]+2{\epsi}^{1/2}.$$
  Secondly,
$\lambda'/r<\lambda'<\frac{1}{(k+1)\hat{T}}$ guarantees \rref{2250} for  $\lambda=\lambda'/r$ because
  $\frac{\ln p}{\lambda}$  is natural.
  At last,  by the definition of $\kappa$ and by the choice of $\hat{T}$ and $\lambda',$ we get
  $$U_{{1}/{\lambda'}}-U_{{r}/{\lambda'}}\leq \kappa({r}/{\lambda'},r)\leq
  \kappa(r/\lambda',p)\leq \epsi/2.$$

  Thus, we obtain
  $$U_{{1}/{\lambda'}}\leq U_{{r}/{\lambda'}}+\epsi/2\leqref{2250}
  V[\bw_{\lambda'/r}]+6{\epsi}^{1/2}+\epsi\leq V[\bw_{\lambda'}]+8{\epsi}^{1/2}+\epsi\leq V[\bw_{\lambda'}]+9{\epsi}^{1/2}$$
for all sufficiently small positive $\lambda'$.
 By arbitrariness  of positive $\epsi$, the proof is complete. \bo

\section{The proof of Proposition \ref{bw_av}.}
\label{B}
 We will continue the notation of the previous section for  values of the game value map~$V$.
 For instance, the symbols
 $\sspk{\int_{0}^h a(t)g(z(t))\,dt+ U_*(z(h))}{\omega}$,
 $\sspk{\int_{0}^h b(t)g(z_1(t))\,dt+ \ssp{c}{z_1(h)}}{z(h')}$
 denote  the values of the game value map 
 for a payoff $\mm{K}\ni z\mapsto \int_{0}^h a(t)g(z(t))\,dt+ U_*(z(h))\in\mm{R}$ at $\omega$
and for  a payoff
 $\mm{K}\ni z_1\mapsto \int_{0}^h b(t)g(z_1(t))\,dt+V[c](z_1(h))\in\mm{R}$ at $z(h')$,
 respectively.

 Also, as in the proof of Proposition \ref{av_bw}, we can  assume that   $\mathfrak{C}$ coincides with the set of all bounded functions $c:\mm{K}\to\mm{R}$ and the game value $V[c]$ is correct for all bounded payoffs $c:\mm{K}\to\mm{R}$ and satisfies \rref{conditions1}--\rref{conditions2}.

{\bf Auxiliary estimates.}\
For each positive $\epsi<1/4$, we can choose  natural  $k\geq 2$ such that $\epsi^{1/2}\leq k\epsi \leq (k+1)\epsi\leq
2{\epsi}^{1/2}.$
Set $p\rav 1+\frac{1}{k}.$
Choose natural $n$ such that $2p^{-n}<\epsi^{1/2}.$
Now,
\begin{eqnarray*}
\frac{\epsi}{p-1}\leq\frac{p\epsi}{p-1}=(k+1)\epsi\leq 2{\epsi}^{1/2},\quad 2(p-1)=\frac{2}{k}\leq 2\epsi^{1/2},\quad
2p^{-n}<\epsi^{1/2}.
\end{eqnarray*}

Set
 $$\kappa(\lambda,p_0)\rav\sup_{p'\in[1,p_0]}
\sup_{\omega\in\Omega}\,\big(\,U_{\lambda}(\omega)-U_{p'\lambda}(\omega)\big)\qquad \forall p_0>1,\lambda>0.$$
By \rref{slowlyl},     there exists positive  $\bar{T}$ such that
 $\kappa(\lambda,p)\leq\epsi/2$ holds for all positive $\lambda<1/\bar{T}$.

  By the definitions of the subsolution and the supersolution,
    there exists positive  $\hat{T}>\bar{T}$ such that,  for all $\omega\in\Omega$, $\lambda\in(0,1/\hat{T})$
      and natural $h>\hat{T},T>h+\hat{T}$,  one has
\begin{eqnarray}
 U_{\lambda}(\omega)&\leq& \ssp{
 \lambda\int_{0}^{h}e^{-\lambda t}g(z(t))\,dt+e^{-\lambda h}U_{\lambda}(z(h))+\epsi/2}{\omega},\label{150131}\\
 V[\av_{T}](\omega)&\geq&\ssp{
     \frac{1}{T}\int_{0}^{h}g(z(t))\,dt+\frac{T-h}{T}V[\av_{T-h}](z(h))-\epsi/3}{\omega}.\label{1501}
\end{eqnarray}
  Note that, for all $\lambda\in (0,1/\hat{T})$, $U_{\lambda}$ is bounded; then,  $e^{-\lambda h_0}U_{\lambda}<\epsi/2$ for all
  sufficiently large natural $h_0$. Now, \rref{150131},\rref{conditions2}, and $g\leq 1$ imply that $U_{\lambda}(\omega)\leq \lambda\int_{0}^{h}e^{-\lambda t}\,dt+\epsi\leq 2$
  for all $\lambda\in (0,1/\hat{T})$.

{\bf The special choice of $T$.}\
 Fix every natural
 $T>kp^{n+1}\hat{T}>k{p^{n+1}}\bar{T}$ such that
 $T(k+1)^{-n-1}$
 is also natural. Then, $T(1-p^{-1})p^{-n}=Tk^n(k+1)^{-n-1}$ is natural as well.
 Set
 $$h\rav T(1-p^{-1})=\frac{T}{k+1},\quad  q\rav p^{-1},\quad
 \lambda\rav\frac{p\ln p}{T(p-1)}=\frac{\ln p}{h}.$$
  In view of $p\ln p>p-1>\ln p$, we also have
 \begin{eqnarray}
 \frac{p}{\lambda}>\frac{p\ln p}{\lambda(p-1)}=T>\frac{1}{\lambda}>\frac{\ln
 p}{\lambda(p-1)}=\frac{T}{p}=\frac{h(k+1)}{p}=\frac{h}{p-1}=kh>kp^{n}\hat{T},\nonumber\\
 p/\lambda>T>1/\lambda>T-h=kh>h>\hat{T}p^{n}>\bar{T}p^{n},\quad e^{-\lambda
 h}=p^{-1}=q=\frac{T-h}{T}=\frac{k}{k+1}.\label{1503_}
\end{eqnarray}
For all $i\in\{0,\dots,n\},$  the numbers $hq^i,(T-h)q^i=khq^i$ are natural; now,
 from $(T-h)q^i>hq^i>\hat{T}$ and $\lambda p^i<1/\bar{T}$  it follows that
  $\lambda p^i,hq^i,Tq^i$ satisfy
the inequalities \rref{150131},\rref{1501}
  for all $i\in\{0,\dots,n\}$.
 Moreover, $U_{p^{i+1}\lambda}\leq 2$ holds  for all $i\in\{0,\dots,n\}$.

{\bf Forward-tracking.}\
By $\lambda<1/\hat{T}<1/\bar{T}$, for all $\omega\in\Omega,$ we have
\begin{eqnarray*}
 {U_{\lambda}(\omega)} &\leqref{150131}&
 \ssp{\int_{0}^{h}\lambda e^{-\lambda t}g(z(t))\,dt+e^{-\lambda h}U_{\lambda}(z(h))+\epsi/2}{\omega}\\
 &\leq&\ssp{\int_{0}^{h}\lambda e^{-\lambda t}g(z(t))\,dt+q U_{p\lambda}(z(h))+\epsi/2+\kappa(\lambda,p)}{\omega}\\
 &\leq&\ssp{\int_{0}^{h}\lambda e^{-\lambda t}g(z(t))\,dt+q U_{p\lambda}(z(h))}{\omega}+\epsi.
\end{eqnarray*}

 Analogously, by $p\lambda<1/\hat{T}<1/\bar{T}$ and $e^{-p\lambda qh}=e^{-\lambda h}=q$, for all $z\in\mm{K}$,
\begin{eqnarray*}
 qU_{p\lambda}(z(h)) &\leq&
 q\ssp{\int_{0}^{qh}p\lambda e^{-p\lambda t}g(z_1(t))\,dt+e^{-p\lambda q h}
 U_{p^2\lambda}(z_1(qh))}{z(h)}+q\epsi/2+q\kappa(p\lambda,p)\\
 &=&\ssp{\int_{0}^{qh}\lambda e^{-p\lambda t} g(z_1(t))\,dt+q^{2} U_{p^2\lambda}(z_1(qh))}{z(h)}+q\epsi.
\end{eqnarray*}
 Substituting it into the relation above, we obtain
 \begin{eqnarray*}
 {U_{\lambda}(\omega)}  &\leq&\ssp{\int_{0}^{h}\!\lambda e^{-\lambda t}g(z(t))\,dt+q U_{p\lambda}(z(h))}{\omega}\!+\epsi\\
 &\leq&\sspll{\int_{0}^{h}\!\lambda e^{-\lambda t} g(z(t))\,dt}
 +\sspk{\int_{0}^{qh}\!\lambda e^{-p\lambda t} g(z_1(t))\,dt+q^{2}
 U_{p^2\lambda}(z_1(qh))}{z(h)}\ssprr{}{\omega}\!\!+\epsi(1+q).
\end{eqnarray*}
 Repeating, by $\lambda p^2<1/\hat{T}<1/\bar{T}$, we get
 \begin{eqnarray*}
 {U_{\lambda}(\omega)}
 &\leq&
 \sspll{\int_{0}^{h}\lambda e^{-\lambda t}g(z(t))\,dt+
 \sspl{\int_{0}^{qh}\lambda e^{-p\lambda t} g(z_1(t))\,dt}}\\
 &\ &
 \sspr{+\sspk{\int_{0}^{q^2h}\lambda e^{-p^2\lambda t} g(z_2(t))\,dt+q^{3}
 U_{p^2\lambda}(z_2(q^2h))}{z_1(qh)}}{z(h)}\ssprr{}{\omega}+\epsi(1+q+q^2).
\end{eqnarray*}
Proceeding in a similar way, in view of $\lambda p^n<1/\hat{T}<1/\bar{T}$, we obtain
 \begin{eqnarray*}
 {U_{\lambda}(\omega)}
 &\leq&\ssplll{\int_{0}^{h}\lambda e^{-\lambda t}g(z(t))\,dt+
 \sspll{\int_{0}^{qh}\lambda e^{-p\lambda t} g(z_1(t))\,dt}}\\
 &\ &+ \sspl{\int_{0}^{q^2h}\lambda e^{-p^2\lambda t} g(z_2(t))\,dt+\dots+\sspl{\int_{0}^{q^n h}\lambda e^{-p^n\lambda t}
 g(z_n(t))\,dt}}\\
 &\ &\ssprrr{\ssprr{\sspr{\sspr{+q^{n+1}    U_{p^{n+1}\lambda}(z_n(q^n h))}{z_{n-1}(q^{n-1}
 h)}\dots}{z_1(qh)}}{z(h)}}{\omega}+\epsi(1+q+q^2+\dots+q^n).
\end{eqnarray*}

{\bf Backtracking.}\
 Note that $U_{p^{n+1}\lambda}\leq 2, g\geq 0$. Then,
 \begin{eqnarray*}
 {U_{\lambda}(\omega)}
 &\leq&\ssplll{\int_{0}^{h}\lambda g(z(t))\,dt+
 \sspll{\int_{0}^{qh}\lambda g(z_1(t))\,dt+ \sspl{\int_{0}^{q^2h}\lambda g(z_2(t))\,dt}}}\\
 &\ &\ssprrr{\ssprr{\sspr{+\dots+\sspk{\int_{0}^{q^n h}\lambda  g(z_n(t))\,dt}{z_{n-1}(q^{n-1}
 h)}\dots}{z_1(qh)}}{z(h)}}{\omega}+2q^{n+1}+\frac{\epsi}{1-q}.
\end{eqnarray*}
 By the choice of $n$   and $p$, we have
  $2q^{n+1}<{\epsi}^{1/2}$ and
  $\frac{\epsi}{1-q}=\frac{p\epsi}{p-1}\leq 2{\epsi}^{1/2}$.
  By \rref{1503_}, $\lambda<p/T$ holds. Thanks to  $V[\av_{q^n (T-h)}]\geq 0$, we obtain
 \begin{eqnarray*}
 {U_{\lambda}(\omega)}
 &\leq&\ssplll{\int_{0}^{h}\lambda g(z(t))\,dt+
 \sspll{\int_{0}^{qh}\lambda g(z_1(t))\,dt+ \sspl{\int_{0}^{q^2h}\lambda g(z_2(t))\,dt}}}\\
 &\ &\ssprrr{\ssprr{\sspr{+\dots+\sspk{\int_{0}^{q^n h}\lambda  g(z_n(t))\,dt}{z_{n-1}(q^{n-1}
 h)}\dots}{z_2(q^2h)}}{z_1(qh)}}{\omega}+3{\epsi}^{1/2}\\
 &\leqref{conditions2}&\ssplll{\frac{p}{T}\int_{0}^{h} g(z(t))\,dt+
 \sspll{\frac{p}{T}\int_{0}^{qh} g(z_1(t))\,dt+\dots+\sspl{\frac{p}{T}\int_{0}^{q^n h} g(z_n(t))\,dt}}}\\
 &\ &\ssprrr{\ssprr{\sspr{
 +q^{n}V[\av_{q^n (T-h)}](z_{n}(q^{n} h))}{z_{n-1}(q^{n-1} h)}\dots}{z(h)}}{\omega}+3{\epsi}^{1/2}\\
 &\ravref{conditions1}&p\ssplll{\frac{1}{T}\int_{0}^{h} g(z(t))\,dt+
 \sspll{\frac{1}{T}\int_{0}^{qh} g(z_1(t))\,dt+\dots+\sspl{\frac{1}{T}\int_{0}^{q^n h} g(z_n(t))\,dt}}}\\
 &\ &\ssprrr{\ssprr{\sspr{
 +q^{n+1}V[\av_{q^n (T-h)}](z_{n}(q^{n} h))}{z_{n-1}(q^{n-1} h)}\dots}{z(h)}}{\omega}+3{\epsi}^{1/2}.
\end{eqnarray*}
  Recall that $q^nT=q^{n-1}(T-h)>\hat{T}$, $q=\frac{q^nT-q^nh}{q^nT}$ (see \rref{1503_}); also, \rref{1501} holds for
  $q^nh,q^nT$. Thus,
 \begin{eqnarray*}
 &\ & \ssp{\frac{1}{T}\int_{0}^{q^n h} g(z_n(t))\,dt
 +q^{n+1}V[\av_{q^n (T-h)}](z_{n}(q^{n} h))}{z_{n-1}(q^{n-1} h)}\\
 &=&q^{n}\ssp{\frac{1}{q^{n}T}\int_{0}^{q^n h} g(z_n(t))\,dt
 +q V[\av_{q^n (T-h)}](z_{n}(q^{n} h))}{z_{n-1}(q^{n-1} h)}\\
 &\leqref{1501}&q^{n}V[\av_{q^n T}]({z_{n-1}(q^{n-1} h)})+q^{n}\epsi/3\\
 &\ravref{1503_}&
 q^{n}V[\av_{q^{n-1} (T-h)}]
 ({z_{n-1}(q^{n-1} h)})+q^{n}\epsi/3.
\end{eqnarray*}
 Then,
 \begin{eqnarray*}
 {U_{\lambda}(\omega)}
 &\leq&p\ssplll{\frac{1}{T}\int_{0}^{h} g(z(t))\,dt+
 \sspll{\frac{1}{T}\int_{0}^{qh} g(z_1(t))\,dt
  +\dots+\sspl{\frac{1}{T}\int_{0}^{q^{n-1} h} g(z_{n-1}(t))\,dt}}}\\
&\ &\ +\ssprrr{\ssprr{\sspr{
 q^{n}V[\av_{q^{n-1} (T-h)}](z_{n-1}(q^{n-1} h))}{z_{n-1}(q^{n-1} h)}\dots}{z(h)}}{\omega}+3{\epsi}^{1/2}+pq^{n}\epsi/3.
\end{eqnarray*}
Proceeding in a similar way, since $q^{i-1}(T-h)\ravref{1503_}q^{i}T$ holds for all $i\in\{0,1,2,\dots,n\}$, we obtain
 \begin{eqnarray*}
 {U_{\lambda}(\omega)}
 &\leq&p\sspk{\frac{1}{T}\int_{0}^{h} g(z(t))\,dt
 +q V[\av_{T-h}](z(h))}{\omega}+3{\epsi}^{1/2}+p(q^{n}+\dots+q)\epsi/3\\
 &\leq&p V[\av_{T}](\omega)+3{\epsi}^{1/2}+p(q^{n}+\dots+q+1)\epsi/3\\
 &\leq&p-1+V[\av_{T}](\omega)+5{\epsi}^{1/2}\\
&\leq& V[\av_{T}](\omega)+7{\epsi}^{1/2}
\end{eqnarray*}

By  \rref{1503_} and $T>\bar{T},$ we get $\kappa(1/T,\frac{p\ln
p}{p-1})\ravref{1503_}\kappa(1/T,T\lambda)\leq\kappa(1/T,p)\leq\epsi/2$. Hence,
\begin{eqnarray}
U_{1/T}\leq U_{\lambda}+\epsi/2\leq V[\av_{T}]+7{\epsi}^{1/2}+\epsi/2\label{2251}
\end{eqnarray}
 for all positive $T> k{p^{n+1}}\hat{T}$ if $T/(k+1)^{n+1}$ is also a natural number.

{\bf The general case.}\
 Consider every positive $T'>k\max\big(p^{n+2}\hat{T},(k+1)^{n+1}\big).$
  Now,  we can choose a positive $r>1$ such that $\frac{T'}{r(k+1)^{n+1}}$  is natural
  and $0\leq T'-\frac{T'}{r}\leq (k+1)^{n+1}$ holds.
  By the choice of $\hat{T}$ and $T',$ we also obtain
  $$1\leq r\leq \bigg(1-\frac{(k+1)^{n+1}}{T'}\bigg)^{-1}\leq
   1+\frac{(k+1)^{n+1}}{T'}< 1+\frac{1}{k}=p.$$
 Now, it follows from
   $T'/r\geq T'/p\geq kp^{n+1}\hat{T}\geq \bar{T}$ that, first,
     by the definition of $\kappa$,
  $$U_{1/T'}-U_{r/T'}\leq \kappa(r/T',r)\leq
  \kappa(r/T',p)\leq \epsi/2,$$
  secondly, \rref{2251} holds for $T=T'/r.$
 At last,
  thanks to  \rref{2012},  we also obtain $$V[\av_{{T'}/{r}}]\leq V[\av_{T'}]+2(r-1)\leq V[\av_{T'}]+2(p-1)\leq
  V[\av_{T'}]+{\epsi}^{1/2}.$$

  Thus,
  $$U_{1/T'}\leq U_{r/T'}+\epsi/2\leqref{2251}
  V[\av_{{T'}/{r}}]+7{\epsi}^{1/2}+\epsi\leq
  V[\av_{T'}]+8{\epsi}^{1/2}+\epsi\leq V[\av_{T'}]+9{\epsi}^{1/2}$$
 for all sufficiently large positive $T'.$
 By arbitrariness  of positive $\epsi$, the proof is complete. \bo
\end{document}